\documentclass[12pt]{amsart} 
\author{Constantin N. Beli}
\title{Decomposability of multivariable
polynomials${}^1$\footnote{$^1$I intend to publish this
paper but first I want to be 100\% sure it is original. If any reader
is aware of any results resembling Theorems 1.3, 4.1, 4.2 or Lemma 3.4
then please let me know. Thank you.}}
\usepackage{amsfonts,amsmath,amssymb}
\def\a{\alpha} \def\b{\beta} \def\c{\gamma} 
 \def\d{\delta} \def\e{\varepsilon} 
  \def\h{\frac} 
\def\j{\infty}  \def\l{\lambda} 
\def\m{\lim}    
\def\p{\partial}

\def\({\overline} \def\){\underline}
\def\<{\cdot} \def\go{\mathfrak}
\def\>{~~~~~~~} \def\#{{\bf
Definition}} \def\be{\begin{equation}}
\def\ee{\end{equation}}

\def\sb{\subset}  \def\sbq{\subseteq}  
\def\ti{\times}   
   
   \def\AA{{\rm a}}
 \def\mo{{\rm mod}~}  
    
\def\p{\go p}    
\def\*{\sharp}  \def\0{} 
 \def\1{^{-1}} \def\rk{{\rm rank}\,} 
 \def\[{\prec} \def\]{\succ} 
\def\bmat{\left(\begin{array}} \def\emat{\end{array}\right)} \def\ev{\equiv}
\def\ap{\cong}

 \def\m2{~(\mo 2)} \def\no{\noindent}
 \def\btm{\begin{thm}}
\def\etm{\end{tm}}
 \def\blem{\begin{lem}}
\def\elem{\end{lem}}

\newtheorem{theorem}{Theorem}[section]
\newtheorem{proposition}[theorem]{Proposition}
\newtheorem{lemma}[theorem]{Lemma}
\newtheorem{definition}{Definition}
\newtheorem{corollary}[theorem]{Corollary}

\newtheorem{bof}[theorem]{}
\newtheorem{teorema}{Theorem}

\def\qed{\mbox{$\Box$}\vspace{\baselineskip}}
\def\pf{$Proof.$} 
\def\bco{\begin{corollary}} \def\eco{\end{corollary}} 
\def\bdf{\begin{definition}} \def\edf{\end{definition}} 
\def\btm{\begin{theorem}} \def\etm{\end{theorem}} 
\def\blm{\begin{lemma}} \def\elm{\end{lemma}} 
\def\bff{\begin{bof}\rm} \def\eff{\end{bof}}
\def\btr{\begin{teorema}} \def\etr{\end{teorema}}
\def\bpr{\begin{proposition}} \def\epr{\end{proposition}}

\def\de{\newcommand} \de\tm[1]{{\no\bf Theorem~#1}}
\def\mb{\mathbb} 
\def\RR{{\mb R}}\def\QQ{{\mb Q}}\def\ZZ{{\mb
Z}}\def\NN{{\mb N}} \def\AA{{\mb A}} \def\PP{{\mb P}}
\def\la{\langle} \def\ra{\rangle}

\de\lm[1]{{\no\bf Lemma~#1}}
\de\df[1]{{\no\bf Definition~#1}} \de\co[1]{{\no\bf Corollary~#1}}

\def\Ff{{\mathcal F}} \def\Vv{{\mathcal V}} 
\DeclareMathOperator\im{Im} 
\DeclareMathOperator\car{char} 
\DeclareMathOperator\conv{conv}

\begin{document}
\maketitle
\begin{quote}
{\bf\footnotesize If $K$ is an algebrically closed field, $n\geq 1$
and $K[X]=K[X_1,\ldots,X_n]$ we determine all finite nonempty sets
$I\sb\NN_0^n$ such that every $P\in K[X]$ of the form $P=\sum_{i\in
I}a_iX^i$ with $a_i\in K^*$ $\forall i\in I$ is decomposable. We also
cosider the problem of finding all sets $I$ such that every
$P=\sum_{i\in I}a_iX^i$ with $a_i\in K^*$ is irreducible and how the
answer to this problem depends on the characteristic.}
\end{quote}
\section{Main theorem, the first implication}

Let $K$ be an algebrically closed field and let $n\geq 1$. Let
$X_1,\ldots,X_n$ be variables. We denote $X=(X_1,\ldots,X_n)$ and
$K[X]=K[X_1,\ldots,X_n]$. If $i=(i_1,\ldots,i_n)$ we denote
$X^i=X_1^{i_1}\cdots X_n^{i_n}$. 

For any $P\in K[X]$ we denote by $I(P)\sb\NN^n$ the support of $P$,
which is the finite set for which $P$ can be written as $P=\sum_{i\in
I(P)}a_iX^i$ with $a_i\in K^*$. (We make the convention that $\NN
=\NN_0$, the set of all non-negative integers.)

Various mathematicians have produced irreducibility criteria in terms
of $I(P)$ only. Most of these results are consequences of the
following principle. If $X_t\nmid P$ $\forall t$ and the Newton
polytope of $P$, $C(P):=\conv I(P)$ cannot be written as a Minkowski
sum $C(P)=C'+C''$, where $C',C''$ are convex polytopes of positive
dimensions with vertices in $\ZZ^n$ then $P$ is
irreducible. However this is not the most general result possible and
the problem of finding {\em all} finite nonempty sets $I\sb\NN^n$ such
that any $P\in K[X]$ with $I(P)=I$ is irreducible is very unlikely to
have a simple solution. In this paper we consider the opposite problem
of determining all sets $I$ such that every $P\in K[X]$ with $I(P)=I$
is decomposable. In our main result, Theorem 1.3, we prove that this
happens in only three trivial cases.

For any finite non-empty $I\sb\NN^n$ we define $V_I=V_I(K)$ and
$Z_I=Z_I(K)$:
$$V_I:=\{ P\in K[X]\mid I(P)=I\}\text{ and }Z_I:=\{ P\in V_I\mid
P\text{ is reducible}\}.$$

We have $V_I\ev (K^*)^{|I|}$, which is an affine variety over $K$. We
will prove that $Z_I$ is a subvariety in $V_I$ and we give neccesary
and sufficient conditions such that $Z_I=V_I$. Obviously if $n=1$,
when every polynomial of degree $\geq 2$ is reducible, we have
$V_I=Z_I$ iff $I\not\sbq\{ 0,1\}$. When $n\geq 2$ the conditions are
much more restrictive. 

We also define:
$$W_I=W_I(K):=\{ P\in K[X]\mid P\neq 0,I(P)\sbq I\}.$$

We have $W_I\ap\AA^{|I|}\setminus\{ 0\}$ so $W_I/K^*\ap\PP^{|I|-1}$. 

For any $P,Q\in K[X]$ we have $I(PQ)\sbq I(P)+I(Q)$. Hence if
$I'+I''\sbq I$ the mapping $(P,Q)\mapsto PQ$ defines a morphism of
quasi-projective varieties $\mu_{I',I''}:W_{I'}\ti W_{I''}\to W_I$. It
induces a morphism of projective varieties
$\(\mu_{I',I''}:W_{I'}/K^*\ti W_{I''}/K^*\to W_I/K^*$.

On $V_I,W_I,W_I/K^*$ we consider the Zariski topology. Then $V_I$ is
open in $W_I$. Also if $I\sbq J$ then $W_I$ is closed in $W_J$ and,
moreover, $W_I$ is the closure of $V_I$ in $W_J$. 

\blm The image of $\mu_{I',I''}:W_{I'}\ti W_{I''}\to W_I$ is closed in
$W_I$. 
\elm
\pf Since $\(\mu_{I',I''}:W_{I'}/K^*\ti W_{I''}/K^*\to W_I/K^*$ is a
morphism of projective varieties $\im\(\mu_{I',I''}$ is closed in
$W_I/K^*$. But the cannonical projection $\pi_I:W_I\to W_I/K^*$ is
continuous and $\im\mu_{I',I''}=\pi_I\1 (\im\(\mu_{I',I''})$ so
$\im\mu_{I',I''}$ is closed in $W_I$.\qed 

\bpr $Z_I$ is closed in $V_I$.
\epr
\pf If $i=(i_1,\ldots,i_n)\in\NN^n$ we denote $\sum i:=i_1+\cdots
+i_n$. For any $d\in\NN$ we denote $I_d=\{ i\in\NN^n\mid\sum i\leq
d\}$. Then $W_{I_d}=\{ P\in K[X]\mid P\neq 0,\deg P\leq d\}$. 

Let $d=\max\{\sum i\mid i\in I\}$. Then $\deg P=d$ for any $P\in
V_I$. We may assume that $d\geq 2$ since otherwise
$Z_I=\emptyset$. Now if $P\in V_I$ we have $P\in Z_I$ iff $P=QR$ for
some $Q,R\in K[X]$ of degree $\leq d-1$, i.e. iff $P$ is in the image
of $\mu_{I_{d-1},I_{d-1}}:W_{I_{d-1}}\ti W_{I_{d-1}}\to
W_{I_{2d-2}}$. Hence $Z_I=V_I\cap\im\mu_{I_{d-1},I_{d-1}}$. But by
Lemma 1.1 $\im\mu_{I_{d-1},I_{d-1}}$ is closed in $W_{I_{2d-2}}$ so
$Z_I$ is closed in $V_I$. \qed

Let $e_1,\ldots,e_n$ be the cannonical basis of $Z^n$,
$e_t=(0,\ldots,1,\ldots,0)$. 

\btm Let $I\sb\NN^n$ be finite and nonempty. Then $Z_i=V_I$ if and
only if one of the following holds.

(i) $I\sb e_t+\NN^n$ and $I\neq\{ e_t\}$ for some $t$.

(ii) There are $i=(i_1,\ldots,i_n)$ and $j=(j_1,\ldots,j_n)$ with
$i_tj_t=0$ $\forall t$ and $\gcd\{i_1,\ldots,i_n,j_1,\ldots,j_n\}
>1$ such that $i,j\in I\sb [i,j]:=\{ (1-\l )i+\l j\mid 0\leq\l\leq
1\}$. 

(iii) $\car K=p$ and $\{ 0\}\neq I\sb p\NN^n$ for some prime $p$.
\etm

{\bf Proof of ``if''.}

Let $P\in V_I$, $P=\sum_{i\in I}a_iX^i$. We consider the three cases.

(i) The condition that $I\sb e_t+\NN^n$ means $X_t\mid P$ and the
condition $I\neq\{ e_t\}$ means that $P$ is not of the form
$aX_t$. Hence $P$ is decomposable.

(ii) Let $\gcd\{i_1,\ldots,i_n,j_1,\ldots,j_n\} =d\geq 2$. Since
$i_tj_t=0$ $\forall t$ we have $d=\gcd\{j_1-i_1,\ldots,j_n-i_n\}$. So
an element $k=i+\l (j-i)\in [i,j]$, with $0\leq\l\leq 1$, will belong
to $\ZZ^n$ iff $\l =\h\a d$ for some integer $0\leq\a\leq d$. So
$I\sbq [i,j]\cap\ZZ^n=\{ i+\h\a d(j-i)\mid 0\leq\a\leq d\}$. Hence
$$P=\sum_{\a =0}^da_{\h{d-\a}di+\h\a dj}X^{\h{d-\a}di+\h\a dj}=\(
P(X^{\h 1di}X^{\h 1dj}),$$
where
$$\( P(Y,Z)=\sum_{\a=0}^da_{\h{d-\a}di+\h\a dj}Y^{d-\a}Z^\a.$$
Note that $a_{\h{d-\a}di+\h\a dj}$ may be $0$ for some $0<\a <d$ but not
for $\a =0,d$ since $i,j\in I=I(P)$. Hence $\( P\in K[Y,Z]$ is a
homogenous polynomial of degree $d\geq 2$ not divisible by $Y$ or
$Z$. Since $K$ is algebrically closed $\( P$ decomposes as $\( P=\(
Q\,\( R$ where $\( Q,\( R\in K[Y,Z]$ are homogenous of degrees
$d',d''\geq 1$ with $d'+d''=d$. So $P=QR$ with $Q=\( Q(X^{\h 1di}X^{\h
1dj})$, $R=\( R(X^{\h 1di}X^{\h 1dj})$. Since $\( Q,\( R$ are not
divisible by $Y$ or $Z$ we have $|I(\( Q)|,|I(\( R)|\geq 2$. Since
$i\neq j$ we have $X^{\h{d'-\a}{d'}i+\h\a{d'}j}\neq
X^{\h{d'-\b}{d'}i+\h\b{d'}j}$ whenever $\a\neq\b$. Hence $|I(Q)|=|I(\(
Q)|\geq 2$ so $Q$ is not a constant. Same for $R$ so $P$ is
decomposable.

(iii) Since $I\sb p\NN^n$ and $\car K=p$ we have $P=Q^p$, where
$Q=\sum_{i\in I}a_i^{\h 1p}X^{\h 1pi}$. But $I\neq\{ 0\}$ so $Q$ is
not a constant so $P$ is decomposable. \qed

\section{Minkowski sums and Newton polytopes}

In this section we give some results on Minkowski sums of convex
polytopes and Newton polytopes which are not esentially new. With
different notations they can be found, say, in [G]. For the sake of
self-containment we include th proofs.

We denote by $\la\cdot,\cdot\ra :\RR^n\ti\RR^n\to\RR$ the usual inner
product. 

If $X\sbq\RR^n$ then we denote by $\conv X$ its convex hull. If
$X,Y\sbq\RR^n$ we denote by $X+Y$ their Minkowski sum, $X+Y:=\{
x+y\mid x\in X,y\in Y\}$. Also if $X\sbq\RR^n$, $y\in\RR^n$ then
$X+y:=X+\{ y\} =\{ x+y\mid x\in X\}$. 

Given $X\sbq\RR$ and $a\in\RR^n$ we denote $\la X,a\ra :=\{\la
x,a\ra\mid x\in X\}$. Assuming that $\min\la X,a\ra$ exists we denote
$$X_a:=\{ x\in X\mid\la x,a\ra =\min\la X,a\ra\}.$$
Throughout this section all subsets of $\RR^n$ will be either finite
sets or polytopes so $\min\la X,a\ra$ will always exists so $X_a$ is
always defined. 

\bff Note that for any $X,Y\sbq\RR^n$ we have $\min\la
X+Y,a\ra=\min\la X,a\ra +\min\la Y,a\ra$ and
$(X+Y)_a=X_a+Y_a$. Also if $x\in X$, $y\in Y$ and $x+y\in(X+Y)_a$ then
$x\in X_a$, $y\in Y_a$.

Similarly if $y\in\RR^n$ then $\min\la
X+y,a\ra=\min\la X,a\ra +\la y,a\ra$ and $(X+y)_a=X_a+y$. 

Any element $x\in\conv X$ can be written as $x=\sum_{t=0}^s\l_tx_t$,
where $x_t\in X$, $\l_t>0$ and $\sum_{t=0}^s\l_t=1$. We have $\la
x,a\ra =\sum_{t=0}^s\l_t\la x_t,a\ra\geq\sum_{t=0}^s\l_t\min\la X,a\ra
=\min\la X,a\ra$, with equality iff $x_t\in X_a$ $\forall t$. Thus
$\min\la\conv X,a\ra =\min\la X,a\ra$ and $(\conv X)_a=\conv X_a$. 
\eff

Let $C$ be a polytope. We have $C=\conv\Vv (C)$, where $\Vv (C)$ is
the set of vertices. We denote by $\Ff (C)$ the set of all faces of
$C$. If $\dim C=m$ then $\Ff (C)=\bigcup_{t=0}^m\Ff_t(C)$, where
$\Ff_t(C)$ is the set of faces of dimension $k$. We have
$\Ff_0(C)=\{\{ v\}\mid v\in\Vv (C)\}$, $\Ff_1(C)$ is the set of edges
and $\Ff_m(C)=\{ C\}$. 

The faces of $C$ are the subsets of $C$ where a linear function
$x\mapsto\la x,a\ra$ reaches its minimum (or maximum). So every
$D\in\Ff (C)$ can be written as $D=C_a$ for some $a\in\RR^n$. 

Let now $C',C''$ be two convex polytopes and let $C=C'+C''$. For
convenience we put $\Vv (C)=\Vv$, $\Ff (C)=\Ff$ and
$\Ff_t(C)=\Ff_t$. Similarly we define $\Vv',\Vv''$, $\Ff',\Ff''$ and
$\Ff_t',\Ff_t''$, corresponding to $C',C''$. 

\bdf For any $D\in\Ff$ we define
$$\Phi (D)=\{ x\in C'\mid\exists y\in C'',~x+y\in D\}$$
$$\Psi (D)=\{ y\in C''\mid\exists x\in C',~x+y\in D\}.$$
Since $\Phi,\Psi$ depend on $C',C''$ we denote them by
$\Phi_{C',C''},\Psi_{C',C''}$. 
\edf

Note that by definition if $D,D'\in\Ff$ and $D'\sbq D$ then $\Phi
(D')\sbq\Phi (D)$ and $\Psi (D')\sbq\Psi (D)$. 

\blm If $D=C_a$ then $\Phi (D)=C_a'$ and $\Psi (D)=C_a''$. 
\elm
\pf If $x\in C_a'$ then for any $y\in C_a''$ we have $x+y\in
C_a'+C_a''=C_a=D$ so $x\in\Phi (D)$. Conversely, if $x\in\Phi
(D)$ then there is $y\in C''$ with $x+y\in D$. We have $x\in C'$,
$y\in C''$ and $x+y\in D=C_a$ so $x\in C_a'$, $y\in C_a''$. Hance
$\Phi (D)=C_a'$. similarly $\Psi (D)=C_a''$. \qed

\bco (i) The mappings $D\mapsto\Phi (D)$, $D\mapsto\Psi (D)$ define
surjective functions $\Phi :\Ff\to\Ff'$ and $\Psi :\Ff\to\Ff''$.

(ii) For any $D\in\Ff$ we have $\Phi (D)+\Psi (D)=D$. For any $x\in
C'$, $y\in C''$ we have $x+y\in D$ iff $x\in\Phi (D)$, $y\in\Psi
(D)$. 
%
\eco
\pf (i) By Lemma 2.3 every $\Phi (D)$ is of the form $C_a'$ for
some $a\in\RR^n$ so it belongs to $\Ff'$. The function $\Phi
:\Ff\to\Ff'$ is surjective because for every $D'\in\Ff'$ there is
some $a\in\RR^n$ such that $D'=C_a'=\phi (C_a)$. Similarly for
$\Psi$. 

(ii) follows from the fact that $C'+C''=C$ and $D=C_a$, $\Phi
(D)=C_a'$ and $\Psi (D)=C_a''$ for some $a\in\RR^n$. \qed

For any $v\in\Vv$ we have $\{ v\}\in\Ff$. From $\{ v\} =\Phi (\{ v\}
)+\Psi (\{ v\} )$ we get that $\Phi (\{ v\} )\in\Ff'$ and $\Psi (\{
v\} )\in\Ff''$ are singeltons. Hence  $\Phi (\{ v\} )=\{ v'\}$ and
$\Psi (\{ v\} )=\{ v''\}$ for some $v'\in\Vv'$, $v''\in\Vv''$. This
leads to the following definition.

\bdf We define $\phi =\phi_{C',C''}:\Vv\to\Vv'$ and $\psi
=\psi_{C',C''}:\Vv\to\Vv''$ by
$$\{\phi (v)\} =\Phi (\{ v\} )\text{ and }\{\psi (v)\} =\Psi (\{
v\} )\text{ }\forall v\in\Vv.$$
\edf 

Note that $\Psi_{C',C''}=\Phi_{C'',C'}$ and
$\psi_{C',C''}=\phi_{C'',C'}$. 

%

\blm (i) $\phi$ and $\psi$ are surjective. 

(ii) If $v\in\Vv$ then $\phi (v)+\psi (v)=v$. Moreover, if $x\in
C'$, $y\in C''$ and $x+y=v$ then $x=\phi (v)$, $y=\psi (v)$. 

(iii) If $D\in\Ff$, $D'=\Phi (D)$ and $D''=\Psi (D)$ then
$\phi_{D',D''}=\phi_{|\Vv (D)}$ and $\psi_{D',D''}=\psi_{|\Vv
(D)}$. Also $\phi (\Vv (D))=\Vv (D')$ and $\psi (\Vv (D))=\Vv (D'')$. 
\elm
\pf (i) Let $v'\in\Vv'$. Then $\{ v'\}\in\Ff'$ so by Corollary 2.3(i)
there is $D\in\Ff$ with $\Phi (D)=\{ v'\}$. Let $v\in\Vv (D)$. Then
$\{ v\}\in\Ff$ and $\{ v\}\sbq D$ so $\Phi (\{ v\} )\sbq\Phi (D)=\{
v'\}$. It follows that $\Phi (\{ v\} )=\{ v'\}$ so $v'=\phi (v)$. Thus
$\phi$ is surjective. Similarly for $\psi$.

(ii) By Corollary 2.3(ii) we have $\{\phi (v)\} +\{\psi (v)\} =\Phi
(\{ v\} )+\Psi (\{ v\} )=\{ v\}$ so $\phi (v)+\psi (v)=v$. For the
second claim note that $x+y=v$ may also be written as $x+y\in\{
v\}$. Since $x\in C'$, $y\in C''$ this implies by Definition 1 that
$x\in\Phi (\{ v\} ) =\{\phi (v)\}$, $y\in\Psi (\{ v\} ) =\{\psi
(v)\}$. 

(iii) We have $D'+D''=D$ so we have the surjections
$\phi':=\phi_{D',D''}:\Vv (D)\to\Vv (D')$ and
$\psi':=\psi_{D',D''}:\Vv (D)\to\Vv (D'')$. If $v\in\Vv (D)$ by (ii)
$\phi'(v)+\psi'(v)=v$. But we have $v\in\Vv$, $\phi'(v)\in C'$ and
$\psi'(v)\in C''$ so, again by (ii), $\phi'(v)=\phi (v)$ and
$\psi'(v)=\psi (v)$. So $\phi'=\phi_{|\Vv (D)}$ and $\psi'=\psi_{|\Vv
(D)}$. Hence $\phi (\Vv (D))=\im\phi'$, which by the surjectivity of
$\phi'$, equals $\Vv (D')$. Similarly $\psi (\Vv (D))=\Vv (D'')$. \qed




\bco (i) If $D\in\Ff$ then $D=\Phi (D)+\Psi (D)$ is the only way of
writing $D$ as $D=X+Y$, with $X,Y$ convex sets, $X\sbq C'$, $Y\sbq
C''$. 

(ii) If $D\in\Ff$, $D'=\Phi (D)$ and $D''=\Psi (D)$ then
$\Phi_{D',D''}=\Phi_{|\Ff (D)}$ and $\Psi_{D',D''}=\Psi_{|\Ff (D)}$. 
\eco
\pf (i) Let $X\sbq C'$, $Y\sbq C''$ be convex sets such that
$X+Y=C$. Then $\forall x\in X$ $\forall y\in Y$ we have $x\in C'$,
$y\in C''$ and $x+y\in X+Y=D$. Thus $x\in\Phi (D)$, $y\in\Psi (D)$ by
Definition 1. So $X\sbq\Phi (D)$, $Y\sbq\Psi (D)$. For the reverse
inclusions let $D'=\Phi (D)$, $D''=\Psi (D)$. Since $D'=\conv\Vv (D')$
$D''=\conv\Vv (D'')$ and $X,Y$ are convex it is enough to prove that
$\Vv (D')\sbq X$, $\Vv (D'')\sbq Y$. But for any $v\in\Vv (D)$ we have
$v\in D=X+Y$ so $v=x+y$ with $x\in X$, $y\in Y$. Since $v\in\Vv$,
$x\in X\sbq C'$ and $y\in Y\sbq C''$ we get $x=\phi (v)$, $y=\psi (v)$
by Lemma 2.4 (ii). Hence $\phi (v)\in X$, $\psi (v)\in Y$ $\forall
v\in\Vv (D)$. It follows that $\Vv (D')=\phi (\Vv (D))\sbq X$ and $\Vv
(D'')=\psi (\Vv (D))\sbq Y$. (See Lemma 2.4(iii).) 

(ii) Since $D=D'+D''$ we have $\Phi':=\Phi_{D',D''}:\Ff (D)\to\Ff
(D')$ and $\Psi':=\Psi_{D',D''}:\Ff (D)\to\Ff (D'')$. If $E\in\Ff (D)$
then $\psi'(E)\sbq D'\sbq C'$ and $\psi'(E)\sbq D''\sbq C''$ are
convex and $E=\Phi'(E)+\Psi'(E)$. Thus $\Phi'(E)=\Phi (E)$ and
$\Psi'(E)=\Psi (E)$ by (i). \qed

\blm (i) For any $[v,w]\in\Ff_1$ there is $0\leq\l_{[v,w]}\leq 1$ such
that $\phi (w)-\phi (v)=\l_{[v,w]}(w-v)$ and $\psi (w)-\psi
(v)=(1-\l_{[v,w]})(w-v)$. 

(ii) For any $[v_1,v_2],\ldots,[v_{N-1},v_N],[v_N,v_1]\in\Ff_1$ we
have\\ $\l_{[v_1,v_2]}(v_2-v_1)+\cdots
+\l_{[v_{N-1},v_N]}(v_N-v_{N-1})+\l_{[v_N,v_1]}(v_1-v_N)=0$. 
\elm
\pf (i) Let $D=[v,w]$. Since $\Vv (D)=\{ v,w\}$ we have by Lemma
2.4(iii) $\Vv (\Phi (D))=\{\phi (v),\phi (w)\}$ so $\Phi (D)=[\phi
(v),\phi (w)]$. (Note that we may have $\phi (v)=\phi (w)$, in which
case $[\phi (v),\phi (w)]$ is a singelton.) Similarly $\Psi (D)=[\psi
(v),\psi (w)]$. By Lemma 2.3(ii) we have $[\phi (v),\phi (w)]+[\psi
(v),\psi (w)]=[v,w]$. In particular, $\phi (w)+\psi (v)\in [v,w]$. But
$\phi (v)+\psi(v)=v$ so $\phi (w)+\psi (v)=v+\phi (w)-\phi (v)$. It
follows that $v+\phi (w)-\phi (v)=v+\l (w-v)$, i.e. $\phi (w)-\phi
(v)=\l (w-v)$ for some $\l =\l_{v,w}$, $0\leq\l\leq 1$. Since also
$(\phi (w)-\phi (v))+(\psi (w)-\psi (v))=w-v$ we get also $\psi
(w)-\psi (v)=(1-\l ) (w-v)$. Obviously $\l_{v,w}=\l_{w,v}$ so in fact
$\l =\l_{[v,w]}$. 

(ii) Follows trivially from (i).\qed

\bpr (i) The polytopes $C',C''$ with $C'+C''=C$ are uniquely
determined up to translations by the coefficients $\l_{[v,w]}$ defined
in Lemma 2.6.

(ii) If $0\leq\l\leq 1$ then $\l_{[v,w]}=\l$ for all $[v,w]\in\Ff_1$
iff $C'=\l C+\b$ and $C''=(1-\l )C-\b$ for some $\b\in\RR^n$.

(iii) $C'$ is a singelton iff all $\l_{[v,w]}$ are $0$. $C''$ is a
singelton iff all $\l_{[v,w]}$ are $1$.
\epr
\pf (i) Let $v_0\in\Vv$ and let $\phi (v_0)=\a$. For any $v\in\Vv$
there is a path of eges going from $v_0$ to $v$, i.e. there are
$v_1,\ldots,v_N=v$ such that $[v_{s-1},v_s]\in\Ff_1$ for $1\leq s\leq
N$. By Lemma 2.6(i) we have $\phi (v)=\a+\sum_{s=1}^N(\phi (v_s)-\phi
(v_{s-1}))=\a+\sum_{s=1}^N\l_{[v_{s-1},v_s]}(v_s-v_{s-1})$. Also $\psi
(v)=v-\phi (v)$. Now $C'=\conv\Vv'=\conv\im\phi$ and
$C''=\conv\Vv''=\conv\im\psi$. Hence $C',C''$ are uniquely determined
by $\a$ and $\l_{[v,w]}$. If $\a$ is replaced by $\a'=\a +\b$ then
$\phi,\psi$ will be replaced by $\phi',\psi'$ with $\phi'(v)=\phi
(v)+\b$ and $\psi'(v)=\psi (v)-\b$. Thus $C'$ and $C''$ will be
replaced by $C_1'=C'+\b$ and $C_1''=C''-\b$.

(ii) If $C'=\l C$, $C''=(1-\l )C$ then $C=C'+C''$ follows from the
convexity of $C$. If $v\in\Vv$ then $v=\l v+(1-\l )v$ and $\l v\in
C',(1-\l )v\in C''$ so $\phi (v)=\l v$, $\psi (v)=(1-\l )v$ by Lemma
2.4(ii). If $[v,w]\in\Ff_1$ then $\l_{[v,w]}(w-v)=\phi (w)-\phi (v)=\l
w-\l v$ so $\l_{[v,w]}=\l$. By (i) all other decompositions $C=C'+C''$
with $\l_{[v,w]}=\l$ $\forall [v,w]\in\Ff_1$ are with $C'=\l C+\b$,
$C''=(1-\l )C-\b$ for som $\b\in\RR^n$.



(iii) follows from (i) with $\l=0$ and $\l =1$. \qed

\bco If $C$ is a triangle then the only decompositions $C=C'+C''$ are
of the form $C'=\l C+\b$, $C''=(1-\l )C-\b$ with $0\leq\l\leq 1$ and
$\b\in\RR^n$. 
\eco
\pf Let $\Vv =\{ u,v,w\}$. By Proposition 2.7 we have to prove that
$\l_{[u,v]}=\l_{[v,w]}=\l_{[w,u]}$. But this follows from
$0=\l_{[u,v]}(v-u)+\l_{[v,w]}(w-v)+\l_{[w.u]}(u-w)=
(\l_{[u,v]}-\l_{[w,u]})(v-u)+(\l_{[v,w]}-\l_{[w,u]})(w-v)$. \qed

\bff{\bf Remarks} The condition (ii) of Lemma 2.6 is equivalent to the
apparently weaker statement that $\l_{[v_1,v_2]}v_2-v_1)+\cdots
+\l_{[v_{N-1},v_N]}(v_N-v_{N-1})+\l_{[v_N,v_1]}(v_1-v_N)=0$ whenever
$v_1,\ldots,v_N$ are the vertices of a 2-dimensional (polygonal) face
of $C$. 

Proposition 2.7 gives only a unicity statement. In fact we also have
existence. Namely the following result holds.
\eff
\btm\footnote{With different notations, this result seems to appear in
[KM, \S4]. It might be older though.}  Let $C$ be a polytope. Let
$\l_{[v,w]}\in [0,1]$ for $[v,w]\in\Ff_1$ with
$\l_{[v_1,v_2]}(v_2-v_1)+\cdots
+\l_{[v_{N-1},v_N]}(v_N-v_{N-1})+\l_{[v_N,v_1]}(v_1-v_N)=0$ whenever
$[v_1,\ldots,v_N]\in\Ff_2$. Then there is a decomposition $C=C'+C''$
such that $\phi :\Vv\to\Vv'$ and $\psi :\Vv\to\Vv''$ satisfy $\phi
(w)-\phi (v)=\l_{[v,w]}(w-v)$ and $\psi (w)-\psi
(v)=(1-\l_{[v,w]})(w-v)$ for all $[v,w]\in\Ff_1$. $C',C''$ are unique
up to translations.
\etm
We don't need this result so we won't prove it here.

\bdf On $\RR^n$ we introduce the partial order relation $\leq$ given
by $(x_1,\ldots,x_n)\leq (y_1,\ldots,y_n)$ if $x_t\leq y_t$ $\forall
t$. For any bounded bellow set $X\sbq\RR^n$ we denote by $\inf X$ its
infinimum. 

We have $\inf X=(i_1,\ldots,i_n)$, where $i_t=\inf\{ x_t\mid
(x_1,\ldots,x_n)\in X\}$.
\edf

Note that if $X$ is finite or a polytope we have in fact $\inf
X=(\min\la X,e_1\ra,\ldots,\min\la X,e_n\ra )$. It follows that $\inf
(X+Y)=\inf X+\inf Y$ and if $y\in\RR^n$ then $\inf (X+y)=\inf
X+y$. Also $\inf \conv X=\inf X$. 

Note that the condition $\inf X\geq 0$ is equivalent to $X\sbq
[0,\infty )^n$. 

\bff{\bf Remark} If $\inf C=0$ and we ask the condition that
$C',C''\sbq [0,\infty )^n$ then in Proposition 2.7 $C',C''$ are
uniquely determined by $\l_{[v,w]}$, not merely uniquely up to
translations. Moreover, $\inf C'=\inf C''=0$. 

Indeed, $C',C''\in [0,\infty )^n$ is equivalent to $\inf C',\inf
C''\geq 0$ and since $\inf C'+\inf C''=\inf C=0$ it is equivalent to
$\inf C'=\inf C''=0$. If $C=C'+C''$ is an arbitrary decomposition
then all the other decompositions corresponding to the same
$\l_{[v,w]}$'s are $C=C_1'+C_1''$ with $C_1'=C'+\b$, $C_1''=C''-\b$ for
some $\b\in\RR^n$. Then $\inf C_1'=\inf C'+\b$ and $\inf C_1''=\inf
C''-\b$. Since $\inf C'+\inf C''=0$ the only value of $\b$ such that
$\inf C_1'=\inf C_1''=0$ is $\b=-\inf C'$. 

Also since $\inf\l C=\inf (1-\l )C=0$ in Proposition 2.7(ii) and
Corollary 2.8 we replace the statement that $C'=\l C+\b$ and
$C''=(1-\l )C-\b$ by $C'=\l C$, $C''=(1-\l )C$.
\eff 

\bff{\bf Integral polytopes.} We say that a polytope is integral if
all its vertices have integer coordinates. We are interested in
decompositions $C=C'+C''$, where both $C',C''$ are integral
polytopes. Then in Proposition 2.7 $\l_{[v,w]}$ should satisfy the
condition $\l_{[v,w]}(w-v)\in\ZZ^n$. This is because $\phi (w)-\phi
(v)\in\ZZ^n$. (We have $\phi (v),\phi (w)\in\Vv'\sbq\ZZ^n$.) The
values of $\l_{[v,w]}\in [0,1]$ satisfying this condition are
precisely $\l_{[v,w]}=\h\a d$, where $d$ is the $\gcd$ of the
coordinates of $w-v$ and $0\leq\a\leq d$. 
\eff

\subsection*{Newton polytopes}
\bdf If $P\in K[X]\setminus\{ 0\}$ we define the Newton polytope of
$P$ as $C(P):=\conv I(P)$.
\edf

\bff Since $C(P)=\conv I(P)$ we have $\Vv (C(P))\sbq I(P)\sb\ZZ^n$ so
$C(P)$ is an integral polytope. We also have $\inf C(P)=\inf I(P)\geq
0$. More precicely $\inf C(P)=\inf I(P)$ is defined by the property
that $X^{-\inf C(P)}P$ is a polynomial not divisible by $X_t$ for any
$t$. 

Also since $\deg X^i=\sum i=\la i,\sum_{t=1}^ne_t\ra$ we have
$$\deg P=\max\la I(P),\sum_{t=1}^ne_t\ra=\max\la
C(P),\sum_{t=1}^ne_t\ra.$$
(Same as for the minimum, we have $\max\la\conv X,a\ra =\max\la
X,a\ra$.) 
\eff

The most important property of the Newton polytope is the following
result first noticed by Ostrowski

\btm{\bf (Ostrowski 1921)} If $P,Q\in K[X]\setminus\{ 0\}$ then
$C(PQ)=C(P)+C(Q)$. 
\etm
\pf Let $C=C(P)+C(Q)$. We have $C=\conv (I(P)+I(Q))$. Since $I(PQ)\sbq
I(P)+I(Q)$ we have $C(PQ)\sbq C$. For the reverse inclusion it is
enough to prove that $\Vv (C)\sbq C(P)+C(Q)$. Let $P=\sum_{i\in
I(P)}a_iX^i$, $Q=\sum_{i\in I(P)}b_iX^i$ with $a_i,b_i\neq 0$. By
Lemma 2.4(ii) every $v\in\Vv (C)$ can be written uniquely as $v=x+y$
with $x\in C(P)$, $y\in C(Q)$, namely $v=\phi (v)+\psi (v)$ with $\phi
(v)\in\Vv (C(P))$, $\psi (v)\in\Vv (C(Q))$. Since $\Vv (C(P))\sbq
I(P)\sbq C(P)$ and $\Vv (C(Q))\sbq I(Q)\sbq C(Q)$ the coefficient of
$X^v$ in $PQ$ is equal to
$$\sum_{i\in I(P),j\in I(Q)\atop i+j=v}a_ib_j=a_{\phi
(v)}b_{\psi(v)}\neq 0.$$ 
Hence $v\in I(PQ)$. \qed 

\section{Proof of the second implication}

We say that a finite nonempty set $I\sb\NN^n$ is good if $Z_I=V_I$. We
want to prove that if $I$ is good then it satisfies one of the
conditions (i)-(iii) of the main theorem. We assume the contrary.

Let $C=\conv I$. Then $C(P)=C$ for any $P\in V_I$. 

\subsection*{The case $\dim C=0$.} In this case $I=C$ is a singleton,
$I=\{ i\}$ and $V_I=\{ aX^i\mid a\in K^*\}$. Hence $Z_I=V_I$ implies
$i\neq 0$ and $i\neq e_t$ $\forall k$. (Otherwise every $P\in V_I$ is
a constant or of the form $aX_t$ so it is irreducible.) Let
$i=(i_1,\ldots,i_n)$ and let $t$ be an index such that $i_t\neq
0$. Then $I=\{ i\}\sb e_t+\NN^n$ and $I\neq\{ e_t\}$ so $I$ satisfies
(i).\qed 

From now on we assume that $\dim C>0$ so $I$ is not a singleton. We
claim that $\inf C=\inf I=0$. Otherwise if $\inf
I=i=(i_1,\ldots,i_n)\neq 0$ then there is some $t$ such that
$i_t>0$. It follows that $I\sb e_t+\NN^n$ and $I\neq\{ e_t\}$ so $I$
satisfies (i). 

\blm If $I$ is good and $\inf C=0$ then there is a decomposition
$C=C'+C''$, where $C',C''$ are integral polytopes of dimension $>0$
and $\inf C'=\inf C''=0$ such that $W_I$ is contained in the image
of $\mu_{C'\cap\ZZ^n,C''\cap\ZZ^n}:W_{C'\cap\ZZ^n}\ti
W_{C''\cap\ZZ^n}\to W_{C\cap\ZZ^n}$. 
\elm
\pf There is a finite number of ways one can write $C=C'+C''$, where
$C',C''$ are integral polytopes and $\inf C'=\inf C''=0$. It is because
by Proposition 2.7 and Remark 2.11 each such pair $(C',C'')$ is
uniquely determined by $\l_{[v,w]}$ with $[v,w]\in\Vv (C)$ and by
2.12 there is a finite number of possible values for each $\l_{[v,w]}$
. Let $(C_1',C_1''),\ldots,(C_N',C_N'')$ be all such pairs for which
also $\dim C',\dim C''>0$. 

Let $P\in V_I=Z_I$. Then $P=QR$, where $Q,R\in K[X]$ are not
constants. Since $\inf C(P)=\inf C=0$ $X_t\nmid P$ $\forall t$ so
$Q,R$ cannot be monomials. Therefore $|I(Q)|.|I(R)|\geq 2$ so $\dim
C(Q),\dim C(R)>0$. Since $C(Q)+C(R)=C(P)=C$ we have
$(C(Q),C(R))=(C_\a',C_\a'')$ for some $\a$. 

We have $I(Q)\sbq C(Q)\cap\ZZ^n=C_\a'\cap\ZZ^n$, $I(R)\sbq
C(R)\cap\ZZ^n=C_\a''\cap\ZZ^n$. Since
$(C_\a'\cap\ZZ^n)+(C_\a'\cap\ZZ^n)\sbq
(C_\a'+C_\a'')\cap\ZZ^n=C\cap\ZZ^n$ the product function
$\mu_{C_\a'\cap\ZZ^n,C_\a''\cap\ZZ^n}:W_{C_\a'\cap\ZZ^n}\ti
W_{C_\a''\cap\ZZ^n}\to W_{C\cap\ZZ^n}$ is defined and $P=QR\in\im
\mu_{C_\a'\cap\ZZ^n,C_\a''\cap\ZZ^n}$. In concluzion
$V_I\sbq\bigcup_{\a =1}^N\im\mu_{C_\a'\cap\ZZ^n,C_\a''\cap\ZZ^n}\sbq
W_{C\cap\ZZ^n}$. Since $\im\mu_{C_\a'\cap\ZZ^n,C_\a''\cap\ZZ^n}$ are
closed in $W_{C\cap\ZZ^n}$ and $V_I$ is a affine variety so it is
connected in the Zariski topology we have
$V_i\sbq\im\mu_{C_\a'\cap\ZZ^n,C_\a''\cap\ZZ^n}$ for some $\a$. By
Lemma 1.1 $\im\mu_{C_\a'\cap\ZZ^n,C_\a''\cap\ZZ^n}$ is closed in
$W_{C\cap\ZZ^n}$ and the closure of $V_I$ in $W_{C\cap\ZZ^n}$ is $W_I$
so $W_i\sbq\im\mu_{C_\a'\cap\ZZ^n,C_\a''\cap\ZZ^n}$. Thus our
statement holds for $(C',C'')=(C_\a',C_\a'')$. \qed

\bff If $C',C''$ are the convex sets from Lemma 3.1 we may use the
notations from \S2. Let $\phi =\phi_{C',C''}:\Vv\to\Vv'$, $\psi
=\psi_{C',C''}:\Vv\to\Vv''$ and let $\l_{[v,w]}$ for $[v,w]\in\Ff_1$
be defined as in Lemma 2.6.

By Proposition 2.7(iii) since neither $C'$ nor $C''$ is a singleton
the numbers $\l_{[v,w]}$, $[v,w]\in\Ff_1$, cannot be all $0$ or all
$1$. 
\eff

\subsection*{The case $\dim C=1$.} In this case we have $C=[i,j]$ for
some $i,j\in\NN^n$ and we have $\{ i,j\}=\Vv (C)\sbq I$ and $I\sb
C=[i,j]$. Let $i=(i_1,\ldots,i_n)$, $j=(j_1,\ldots,j_n)$. Since
$0=\inf C=\inf\{ i,j\}$ ($C=\conv\{ i,j\}$) we have $i_tj_t=0$
$\forall t$. 

By 2.12 we have $\l_{[i,j]}=\h\a d$ with $0\leq\a\leq d$, where
$d=\gcd\{ j_1-i_1,\ldots,j_n-i_n\}$ Since $[i,j]$ is the only edge of
$C$ we have by 3.2 $0<\l_{[i,j]}<1$ so $0<\a <d$, which implies that
$d>1$. But $i_tj_t=0$ $\forall t$ so $gcd\{
i_1,\ldots,i_n,j_1,\ldots,j_n\} =\gcd\{ j_1-i_1,\ldots j_n-i_n\}
=d>1$.

Since also $i,j\in I\sbq [i,j]$ $I$ satisfies (ii). \qed

Before proving the case $\dim C\geq 2$ we need two more Lemmas.

\blm (i) With the notations from 3.2, we have $\l_{[i,j]}=\l$ $\forall
[i,j]\in\Ff_1$ for some $0<\l <1$. Hence $C'=\l C$, $C''=(1-\l )C$. 

(ii) If $[i,j]\in\Ff_1$ and $\{ i,j\}\sbq J\sbq I$ then $J-\inf J$ is
a good set.
\elm
\pf First note that if $i\in\Vv$ and $i\in J\sbq I$ then $i$ is also a
vertex of $\conv J$. Indeed, we have $\{ i\}\in\Ff$ so there is
$a\in\RR^n$ with $\{ i\} =C_a$ i.e. $\la i,a\ra <\la x,a\ra$ for any
$x\in C\setminus\{ i\}$. Since $i\in\conv J\sbq\conv I=C$ we have $\la
i,a\ra <\la x,a\ra$ for any $x\in\conv J\setminus\{ i\}$ so $\{
i\}=(\conv J)_a$ so $i$ is a vertex of $\conv J$.

Assume now that $P\in V_J$. Then $P\in
W_I\sbq\im\mu_{C'\cup\ZZ^n,C''\cup\ZZ^n}$ so there are $Q\in
W_{C'\cup\ZZ^n}$, $R\in W_{C''\cup\ZZ^n}$ such that $P=QR$. If
$C_1=C(P)=\conv J$, $C_1'=C(Q)$ and $C_1''=C(R)$ then $C_1=C_1'+C_1''$
and we denote by $\phi'=\phi_{C_1',C_1''}:\Vv (C_1)\to\Vv (C_1')$,
$\psi'=\psi_{C_1',C_1''}:\Vv (C_1)\to\Vv (C_1'')$ and $\l_{[v,w]}'$ for
$[v,w]\in\Ff_1 (C_1)$ the $\phi,\psi$ and $\l_{[v,w]}$ corresponding
to $C_1,C_1',C_1''$. Now $I(P)\sbq I$ so $C_1\sbq C$, $I(Q)\sbq C'$ so
$C_1'\sbq C'$ and $I(R)\sbq C''$ so $C_1''\sbq C''$. We have $i\in\Vv
(C_1)$ so the only way of writing $i=x+y$ with $x\in C_1'$, $y\in
C_1''$ is $i=\phi'(x)+\psi'(y)$. On the other hand $i\in\Vv$ so the
only way of writing $i=x+y$ with $x\in C'$, $y\in C''$ is
$i=\phi(x)+\psi(y)$. Thus $\phi'(i)=\phi (i)$ and $\psi'(i)=\psi
(i)$. 

After this introduction we start the proof of our lemma.

(i) Consider first two edges of $C$ that meet in the same vertex,
$[i,j],[j,k]\in\Ff_1$. Take $J=\{ i,j,k\}$ and construct
$C_1,C_1',C_1'',\phi',\psi'$ and $\l'_{[v,w]}$ for
$[v,w]\in\Ff_1(C_1)$ as above. Now $C_1=\conv J$ is the triangle
$[i,j,k]$ so by Proposition 2.7(ii) and Corollary 2.8 we have
$\l_{[i,j]}'=\l_{[j,k]}'=\l_{[k,i]}'$. But $i,j,k\in\Vv$ so
$\phi'(i)=\phi (i)$, $\phi'(j)=\phi (j)$, $\phi'(k)=\phi (k)$. Thus
$\l_{[i,j]}'(j-i)=\phi'(j)-\phi'(i)=\phi (j)-\phi
(i)=\l_{[i,j]}(j-i)$ so $\l_{[i,j]}'=\l_{[i,j]}$. Similarly
$\l_{[j,k]}'=\l_{[j,k]}$ so $\l_{[i,j]}=\l_{[j,k]}$.

Let now $[i,j],[i',j']\in\Ff_1$ be arbitrary. Then there is a path
between $j$ and $i'$ along the edges of $C$,
$i_1=j,i_2,\ldots,i_N=i'$. Then $\l_{[i,i_1]}=\l_{[i_1,i_2]}=\ldots
=\l_{[i_{N-1},i_N]}=\l_{[i_N,j']}$ so $\l_{[i,j]}=\l_{[i',j']}$.

Let $\l$ be the common value of all $\l_{[i,j]}$. By 3.2 $\l\neq 0,1$
so $0<\l <1$. 

(ii) Note that $\inf (J-\inf J)=\inf J-\inf J=0$ so $J-\inf
J\sb\NN^n$. Let $\( P\in V_{J-\inf J}$ and let $P=X^{\inf J}\(
P$. Then $P\in V_J$. We then decompose $P=QR$ as above. Since
$i,j\in\Vv$ and $i,j\in J$ we have $i,j\in\Vv (C_1)$ and
$\phi'(i)=\phi (i)$, $\psi'(i)=\psi (i)$, $\phi'(j)=\phi (j)$,
$\psi'(j)=\psi (j)$. By (i) $\l_{[i,j]}=\l\neq 0,1$ so $\phi (j)-\phi
(i)=\l (j-i)\neq 0$ and $\psi (j)-\psi (i)=(1-\l )(j-i)\neq 0$. Hence
$\phi'(i)\neq\phi'(j)$ and $\psi'(i)\neq\psi'(j)$. Since
$\phi'(i),\phi'(j)\in C_1'$ and $\psi'(i),\psi'(j)\in C_1'$ neither
$C_1'=C(Q)$, nor $C_1''=C(R)$ is a singleton so neither $Q$, nor $R$
is a monomial. Since $X^{\inf J}\( P=P$ is the product of two
polynomials, neither of which is a monomial, we have that $\( P$ is
decomposable so $\( P\in Z_{J-\inf J}$. Hence $Z_{J-\inf J}=V_{J-\inf
J}$, i.e. $J-\inf J$ is good. \qed

\blm Let $i,j,k\in\NN^n$ be noncolinear and let $I=\{ i,j,k\}$. Assume
that $\inf I=0$. Then $Z_I=V_I$ if $\car K=p$ and $I\sb p\ZZ^n$ for
some prime $p$ and $Z_I=\emptyset$ otherwise. 
\elm
\pf The first statement follows trivially from the first section as
$I$ satisfies condition (iii) of the main theorem. 

For the second statement let $P=a_iX^i+a_jX^j+a_kX^k$ where
$a_i,a_j,a_k\in K^*$. Assume that $P\in Z_I$, i.e. $P$ decomposes. We
write $P=QR$, where $Q,R\in K[X]$ are not constants and $Q$ is
irreducible. If $C=C(P)$, $C'=C(Q)$ and $C''=C(R)$ then
$C=C'+C''$. Since $C$ is a triangle with $\inf C=\inf I=0$ and
$C',C''\sbq [0,\j )^n$ we have by Corollary 2.8 and 2.11 $C'=\l C$,
$C''=(1-\l )C$ with $\l\in [0,1]$. Since $Q,R$ are not constants we
have $C',C''\neq\{ 0\}$ so $\l\neq 0,1$. 

Let $Q=\sum_{h\in I(Q)}b_hX^h$ with $b_h\in K^*$. Now $C(Q)=C'=\l C$
is the triangle with vertices at $\l i,\l j,\l k$ so $\l i,\l j,\l
k\in I(Q)\sb\NN ^n$. It follows that $\l\in\QQ$ and if we write $\l
=\h\a d$ with $(\a,d)=1$ then $i,j,k\in d\NN^n$, i.e. $I\sb
d\NN^n$. By hypothesis if $\car K=p$ then $I\not\sb p\NN^n$ so $p\nmid
d$. 

Before going further we prove that if $s_1,t_1,s_2,t_2\in\ZZ$ with
$s_1t_2-s_2t_1\neq 0$ and $\e_1,\e_2\in K^*$ then the system
$x^{s_1}y^{t_1}=\e_1$, $x^{s_2}y^{t_2}=\e_2$ has a solution with
$x,y\in K^*$. If $t_2\neq 0$ we note that for any $q\in\ZZ$ our system
is equivalent to the system $x^{s_2}y^{t_2}=\e_2$,
$x^{s_1}y^{t_1}(x^{s_2}y^{t_2})^{-q}=\e_1\e_2^{-q}$. The second
equation may be writen as $x^{s_3}y^{t_3}=\e_3$, where
$\e_3=\e_1\e_2^{-q}$, $s_3=s_1-qs_2$ and $t_3=t_1-qt_2$. We choose $q$
such that $|t_3|<|t_2|$. Note that $s_2t_3-s_3t_2=-(s_1t_2-s_2t_1)\neq
0$. We repeat the procedure and obtain new equivalent systems
$x^{s_l}y^{t_l}=\e_l$, $x^{s_{l+1}}y^{t_{l+1}}=\e_{l+1}$ with
$s_lt_{l+1}-s_{l+1}t_l\neq 0$ and $|t_2|>\ldots >|t_{l+1}|$ until we
get an index $l$ with $t_{l+1}=0$. We have $0\neq
s_lt_{l+1}-s_{l+1}t_l=-s_{l+1}t_l$ so $s_{l+1},t_l\neq 0$ and our
system writes as $x^{s_l}y^{t_l}=\e_l$, $x^{s_{l+1}}=\e_{l+1}$, which
obviously has solutions. (Take $x$ with $x^{s_{l+1}}=\e_{l+1}$ and
then take $y$ with $y^{t_l}=\e_lx^{-s_l}$.)



Let now $q=(q_1,\ldots,q_n),r=(r_1,\ldots,r_n)\in\ZZ^n$ be linearly
independent and let $\e,\eta\in K^*$. We claim that there is
$x=(x_1,\ldots,x_n)\in (K^*)^n$ such that $x^q=\e$, $x^r=\eta$. 

Since the $2\ti n$ matrix $q\choose r$ has rank $2$ there are $1\leq\a
<\b\leq n$ such that $q_\a r_\b -q_\b r_\a\neq 0$. Then as seen above
there are $x,y\in K^*$ such that $x^{q_\a}y^{q_\b}=\e$ and
$x^{r_\a}y^{r_\b}=\eta$. Then we simply take $x=(x_1,\ldots,x_n)$ with
$x_\a =x$, $x_\b =y$ and $x_t=1$ for $t\neq\a,\b$ and we have
$x^q=\e$, $x^r=\eta$. 

We now continue our proof. Since $\h 1d(j-i)$ and $\h 1d(k-i)$ are
linearly independent for any $\e,\eta\in\mu_d$ there is $\c (\e,\eta
)=(\c_1(\e,\eta ),\ldots\c_n(\e,\eta ))\in (K^*)^n$ such that $\c
(\e,\eta )^{\h 1d(j-i)}=\e$, $\c (\e,\eta )^{\h 1d(k-i)}=\eta$. We
define $P_{\e,\eta}:=P(\c_1(\e,\eta )X_1,\ldots,\c_n(\e,\eta )X_n)$
and $Q_{\e,\eta}:=Q(\c_1(\e,\eta )X_1,\ldots,\c_n(\e,\eta
)X_n)$. Since $Q\mid P$ we have $Q_{\e,\eta}\mid P_{\e,\eta}$ and since
$Q$ is irreducible so is $Q_{\e,\eta}$.

For any $h\in\ZZ^n$ we have $(\c_1(\e,\eta )X_1,\ldots,\c_n(\e,\eta
)X_n)^h=\c (\e,\eta )^hX^h$ so $P_{\e,\eta}=a_i\c (\e,\eta)^iX^i+a_j\c
(\e,\eta)^iX^j+a_l\c (\e,\eta)^kX^k$ and\\ $Q_{\e,\eta}=\sum_{h\in
I(Q)}b_h\c (\e,\eta )^hX^h$.  

Note that $P_{\e,\eta}=\c (\e,\eta )^i(a_iX^i+a_j\c (\e,\eta
)^{j-i}X^j+a_k\c (\e,\eta )^{k-i}X^k)$. But $\c (\e,\eta
)^{j-i}=\e^d=1$ and $\c (\e,\eta )^{k-i}=\eta^d=1$ so $P_{\e,\eta}=\c
(\e,\eta )^iP\sim P$. (Here by $\sim$ we mean that the two polynomials
are associates in $K[X]$.) Since $Q_{\e,\eta}\mid P_{\e,\eta}$ we have
$Q_{\e,\eta}\mid P$ $\forall\e,\eta\in\mu_d$. Assume now that
$Q_{\e,\eta}\sim Q_{\e',\eta'}$ so $Q_{\e',\eta'}=tQ_{\e,\eta}$ for
some $t\in K^*$. We have $\l i,\l j,\l k\in I(Q)$ and by considering
the coefficients of $X^{\l i},X^{\l j},X^{\l k}$ we get $t=\h{b_i\c
(\e',\eta')^{\l i}}{b_i\c (\e,\eta )^{\l i}}=\h{b_j\c (\e',\eta')^{\l
j}}{b_j\c (\e,\eta )^{\l j}}=\h{b_l\c (\e',\eta')^{\l
k}}{b_l\c (\e,\eta )^{\l k}}$. It follows that $\c (\e,\eta )^{\l
(j-i)}=\c (\e',\eta')^{\l (j-i)}$ and $\c (\e,\eta )^{\l
(k-i)}=\c (\e',\eta')^{\l (k-i)}$ so $\e^\a =\e'^\a$ and $\eta^\a
=\eta'^\a$. (Recall, $\l =\h\a d$.) But $\e,\e',\eta,\eta'\in\mu_d$
and $(\a,d)=1$ so $\e =\e'$ and $\eta =\eta'$. Since $Q_{\e,\eta}\mid
P$ $\forall\e,\eta\in\mu_d$ and $Q_{\e,\eta}\not\sim Q_{\e',\eta'}$ if
$(\e,\eta)\neq (\e',\eta')$ we have
$\prod_{\e,\eta\in\mu_d}Q_{\e,\eta}\mid P$. It follows that $\deg
P\geq\deg\prod_{\e,\eta\in\mu_d}Q_{\e,\eta}=|\mu_d|^2\deg Q$. But
$C(P)=C$ and $C(Q)=\l C$ so by 2.13 $\deg P=\max\la
C,\sum_{t=1}^ne_t\ra$ and $\deg Q=\max\la\l C,\sum_{t=1}^ne_t\ra
=\l\deg P\geq\h 1d\deg P$. Since $\car K\nmid d$ we also have
$|\mu_d|=|\mu_d(K)|=d$ so $|\mu_d|^2\deg Q\geq d^2\cdot\h 1d\deg
P=d\deg P>\deg P$. Contradiction. Hence $P$ is irreductible $\forall
P\in V_I$, i.e. $Z_I=\emptyset$. \qed

\subsection*{The case $\dim C\geq 2$.} We will prove that $I$
satisfies condition (iii) of the main theorem, i.e. that $\car K=p$
and $I\sb p\NN^n$ for some prime $p$. In fact we only have to prove
that all elements of $I$ are $\ev\mod p\ZZ^n$. Indeed, this implies
that for any index $t$ the $t$-th coordinates of all elements of $I$
are $\ev\mod p$. But $\inf I=0$ so at least one of these $t$-th
coordinates is $0$. Hence all $t$-th coodinates are divisible by
$p$. So $I\sbq p\NN^n$.

Let $i\in\Vv$. We will show that for any $k\in I$, $k\neq i$ we have
$\car K=p$ and $i\ev k\pmod{p\ZZ^n}$ for some prime $p$. Let
$[i,j],[i,j']\in\Ff_1$ be two different edges emerging from $i$. Since
$k\notin [i,j]\cap [i,j']=\{ i\}$ we may assume that, say, $k\notin
[i,j]$. Since $k\in C$ this implies that $k$ is not colinear with
$i,j$. We use Lemma 3.3(ii) with $J=\{ i,j,k\}$. If $h=\inf J$ then
$J-h=\{ i-h,j-h,k-h\}$ is a good set, i.e. $Z_{J-h}=V_{J-h}$ and we
have $\inf (J-h)=h-h=0$. Since $i-h,j-h,k-h$ are not colinear we use
Lemma 3.3 and we have $\car K=p$ and $J-h\sb p\NN^n$ for some prime
$p$. Since $i-h,k-h\in p\NN^n$ we get $i\ev k\pmod{p\ZZ^n}$. \qed

\section{Idecomposability and characteristic}

From Theorem 1.3 we see that if $I\sb\NN^n$ is a finite nonempty set
then $Z_I=V_I$ holds either for every field $K$, if (i) or (ii) holds,
or only for fields with $\car K\in S$, where $S$ is a finite, posibly
empty set of primes depending on $I$, otherwise. More precisely, if
$I$ doesn't satisfy (i) or (ii) then $S=\emptyset$ if $I=\{ 0\}$ and
$S=\{ p\mid I\sb p\NN^n\}$, otherwise.

As we mentioned in \S1, given a finite nonempty set $I\sb\NN^n$, the
question if $Z_I=\emptyset$ doesn't seem to have a simple answer and
we will not attempt to solve this problem. We will prove however that,
same as for the question if $Z_I=V_I$, the answer depends on the
characteristic we show the nature of this dependence. 

\btm Let $I\sb\NN^n$ be finite and nonempty. Then there is a finite,
possibly empty, set of primes $S$ satisfying one of the following.

(a) $Z_I=\emptyset$ iff $\car K\in S$.

(b) $Z_I=\emptyset$ iff $\car K\not\in S$.
\etm
\pf We use the notation from the proof of Proposition 1.2. The proof
relies on Hilbert's Nullstellensatz. For every $P\in V_I$ we have
$\deg P=d:=\max\{\Sigma i\mid i\in I\}$. Such a polynomial is
decomposable iff there are $Q,R$ of degree $\leq d-1$ with
$QR=P$. Polynomials $Q,R$ may be written as $Q=\sum_{i\in
I_{d-1}}a_iX^i$, $R=\sum_{j\in I_{d-1}}b_jX^j$ and we have
$P=\sum_{k\in I_{2d-2}}c_kX^k$ where $c_k=\sum_{i,j\in
I_{d-1},i+j=k}a_ib_j$. The condition that $P\in V_I$, i.e. $I(P)=I$
means $c_k\neq 0$ if $k\in I$ and $c_k=0$ if $k\in I_{2d-2}\setminus
I$. Hence $Z_I\neq 0$, i.e. $V_I$ contains decomposable polynomials
iff there are $a=(a_i)_{i\in I_{d-1}},b=(b_i)_{i\in I_{d-1}}\in
K^{I_{d-1}}$ such that $c_k\neq 0$ if $k\in I$ and $c_k=0$ if $k\in
I_{2d-2}\setminus I$. 

Let $A=(A_i)_{i\in I_{d-1}}$ and $B=(B_j)_{j\in I_{d-1}}$ be
multivariables. For any $k\in I_{2d-2}$ we consider the polynomial
$f_k\in K[A,B]$, $f_k=\sum_{i,j\in I_{d-1},i+j=k}A_iB_j$. Then $c_k$
above can be written as $c_k=f_k(a,b)$. Hence $Z_I\neq\emptyset$ iff
there are $a,b\in K^{I_{d-1}}$ such that $f_k(a,b)\neq 0$ if $k\in I$,
$f_k(a,b)=0$ if $k\in I_{2d-2}\setminus I$. Equivalently,
$Z_I=\emptyset$ iff for any $a,b\in K^{I_{d-1}}$ such that
$f_k(a,b)=0$ $\forall k\in I_{2d-2}\setminus I$ we have $\prod_{k\in
I}f_k(a,b)=0$. By Nullstellensatz there is some $N\in\NN$ and there
are $g_k\in K[A,B]$ for $k\in I_{2d-2}\setminus I$ such that
$$(\prod_{k\in I}f_k)^N=\sum_{k\in I_{2d-2}\setminus I}f_kg_k.$$

Moreover, by the effective Nullstellensatz $N$ and an upper bound $D$
for the degrees of $g_k$ can be found in terms of $I$ alone. Then the
equation above is equivalent to a linear system where the unknowns are
the coefficients of the $g_k$'s. Since all $f_k$'s have integer
coefficients the augmented matrix of the system has integer
entries. Let $U\in M_{q,r}(\ZZ )$ and $\( U\in M_{q,r+1}(\ZZ )$ be the
matrix and the augmented matrix and let $\d =\rk\( U-\rk U\in\{
0,1\}$. For a given field $K$ the matrix and the augmented matrix are
$U(K)$ and $\( U(K)$, the images of $U$ and $\( U$ in $M_{q,r}(K)$ and
$M_{q,r+1}(K)$. Let $\d (K)=\rk\( U(K)-\rk U(K)$. If $\d (U)=0$
then the system is compatible and so $Z_I=\emptyset$, while if $\d
(K)=1$ it is incompatible so $Z_I\neq\emptyset$. But $U(K)$ and $\(
U(K)$ and so $\delta (K)$ only depend on the characteristic $\kappa$
of $K$ so we will denote them by $U_\kappa,\( U_\kappa,\d_\kappa$. If
$\kappa =0$ then $U_0=U$ and $\( U_0=\( U$ so $\d_0=\d$. If $\kappa
=p>0$ then $U_p$ and $\( U_p$ are the images of $U$ and $\( U$ in
$M_{q,r}(\ZZ/p\ZZ )$ and $M_{q,r+1}(\ZZ/p\ZZ )$ and $\d_p=\rk\(
U_p-\rk U_p$. For almost all $p$ we have $\rk U_p=\rk U$ and
$\rk\( U_p=\rk\( U$ so $\d_p=\d$. We denote by $S$ the finite set
of all primes $p$ with $\d_p\neq\d$. Then we have either $\d_\kappa
=0$ if $\kappa\in S$ and $\d_\kappa =1$ if $\kappa\not\in S$, so (a)
holds, or $\d_\kappa =1$ if $\kappa\in S$ and $\d_\kappa =0$ if
$\kappa\not\in S$, so (b) holds. \qed

\btm For any finite set of primes $S$ each of the cases (a) and (b)
of Theorem 4.1 can be obtained for some $I$.
\etm
\pf For (b) we use Lemma 3.4. Let $d=\prod_{p\in S}p$. We define $I=\{
(0,0),(d,0),(0,d)\}\sb\NN^2$. Since the three points are not colinear
and $\inf I=0$ Lemma 3.4 applies. The only primes $p$ with $I\sb
p\NN^2$ are precisely those in $S$. So by Lemma 3.4
$Z_I=V_I\neq\emptyset$ if $p\in S$ and $Z_I=\emptyset$
otherwise. Hence (b) holds. 

If we take, say, $I=\{ 1,2\}$ then $Z_I=V_I\neq\emptyset$ regardless
of the characteristic of $K$ so (a) holds with $S=\emptyset$. However
constructing sets $I$ satisfying (a) for $S\neq\emptyset$ is more
difficult and requires some preliminary results.
\vskip 3mm

If $I\sb\NN^n$ is finite and nonempty and $K$ is an algebraically
closed field we denote by $A(I,K)$ the statement
$Z_I(K)=\emptyset$. If $I',I'',I\sb\NN^n$ are finite and nonempty
denote by $B(I',I'',I,K)$ the statement: There are no $P,Q\in K[X]$
with $I(P)=I'$, $I(Q)=I''$, $I(PQ)=I$. 

\blm If $I',I'',I\sb\NN^n$ are finite and nonempty and $K$ is an
algebraically closed field then $B(I',I'',I,K)$ is equivalent to
$A(J,K)$, where $J\sb\NN^{n+2}$,
$$J=(I\ti\{ (0,0)\})\cup (I'\ti\{ (0,1)\} )\cup(I''\ti\{
(1,0)\})\cup\{ (0,\ldots,0,1,1)\}.$$ 
\elm
\pf We introduce two more variables, $X_{n+1}=Y$ and $X_{n+2}=Z$.

If $B(I',I'',I,K)$ fails then let $P,Q\in K[X]$ with $I(P)=I'$,
$I(Q)=I''$, $I(PQ)=I$. Then $I((P+Y)(Q+Z))=I(PQ+PZ+QY+YZ)=(I(PQ)\ti\{
(0,0)\} )\cup(I(P)\ti\{ (0,1)\} )\cup(I(Q)\ti\{ (1,0)\} )\cup\{
(0,\ldots,0,1,1)\} =J$. Hence $(P+Y)(Q+Z)\in Z_J(K)$ so $A(J,K)$
fails. 

Conversely, if $A(J,K)$ fails let $F\in Z_J(K)$. We have $I(F)=J$ so
$F=R+PZ+QY+YZ$ for some $P,Q,R\in K[X]$ with $I(R)=I$, $I(P)=I'$,
$I(Q)=I''$. Now $F$ is decomposable as $F=GH$, with $G,H\in
K[X,Y,Z]\setminus K$. By considering the degrees in $Y$ and $Z$ and
the coefficient of $YZ$ one shows that, up to permutation and
multiplication by constants, we have $G=P'+Y$, $H=Q'+Z$ with $P',Q'\in
K[X]$. From $R+PZ+QY+YZ=(P'+Y)(Q'+Z)$ one gets $P=P'$, $Q=Q'$ and
$R=P'Q'$ so $R=PQ$. Since $I(P)=I(P')=I'$, $I(Q)=I(Q')=I''$ and
$I(PQ)=I(R)=I$ $B(I',I'',I,K)$ fails. \qed

\blm Let $d\geq 2$ be an integer. We define $I',I'',I\sb\NN^3$ as
follows.
$$I'=(\{ 0,1\}\ti\{ 0,1\}\ti\{ 0\} )\cup (\{ 0\}\ti\{ 0,1\}\ti\{ 1\}
),$$ 
$$I''=\{ 0,\ldots,d\}\ti\{ 0,1\}\ti\{ 0\}$$
and
$$I=(I'+I'')\setminus (\{ 1,\ldots,d\}\ti\{ 0,2\}\ti\{ 0\}
)\setminus\{ (0,1,1),(d,1,1)\}.$$

Then $B(I',I'',I,K)$ holds if and only if $\car K\mid d$.
\elm
\pf We have 
$$I'+I''=(\{ 0,\ldots,d+1\}\ti\{ 0,1,2\}\ti\{ 0\} )\cup (\{
0,\ldots,d\}\ti\{ 0,1,2\}\ti\{ 1\} ).$$ 

We denote $(X_1,X_2,X_3)=(X,Y,Z)$.

If $P,Q\in K[X,Y,Z]$ with $I(P)=I'$, $I(Q)=I''$ then $P=P_0+P_1Z$,
where 
$$P_0=\sum_{i=0}^1\sum_{j=0}^1a_{i,j}X^iY^j\text{ and
}P_1=a_0'+a_1'Y$$
and
$$Q=\sum_{i=0}^d\sum_{j=0}^1b_{i,j}X^iY^j,$$
where all $a_{i,j}$, $a'_j$, $b_{i,j}$ belong to $K^*$.

Since $I(PQ)\sbq I'+I''$ we have $PQ=R_0+R_1Z$, where
$$R_0=\sum_{i=0}^{d+1}\sum_{j=0}^2c_{i,j}X^iY^j\text{ and
}R_1=\sum_{i=0}^d\sum_{j=0}^2c_{i,j}'X^iY^j.$$

Since $(P_0+P_1Z)Q=R_0+R_1Z$ we have $P_0Q=R_0$ so
$$(a_{0,0}+a_{1,0}X)\sum_{i=0}^db_{i,0}X^i=
\sum_{i=0}^{d+1}c_{i,0}X^i,~
(a_{0,1}+a_{1,1}X)\sum_{i=0}^db_{i,1}X^i=
\sum_{i=0}^{d+1}c_{i,2}X^i$$
$$\text{and }(a_{0,0}+a_{1,0}X)\sum_{i=0}^db_{i,1}X^i+
(a_{0,1}+a_{1,1}X)\sum_{i=0}^db_{i,0}X^i=
\sum_{i=0}^{d+1}c_{i,1}X^i$$
and $P_1Q=R_1$ so
$$(a_0'+a_1'Y)(b_{i,0}+b_{i,1}Y)=c_{i,0}'+c_{i,1}'Y+c_{i,2}'Y^2\text{
for }0\leq i\leq d.$$

Suppose now that $\car K=p\mid d$ and $B(I',I'',I,K)$ fails so there
are $P,Q\in K[X,Y,Z]$ with $I(P)=I'$, $I(Q)=I''$ and $I(PQ)=I$. We
have $d=\p^td'$ with $t>0$ and $p\nmid d'$. Using the above notations
the condition that $I(PQ)=I$ means $c_{1,j}=\ldots =c_{d,j}=0$ for
$j=0,2$, $c_{0,1}'=c_{d,1}'=0$ and all other $c_{i,j}$ and $c_{i,j}'$
are $\neq 0$.

We have
$(a_{0,0}+a_{1,0}X)\sum_{i=0}^db_{i,0}X^i=\sum_{i=0}^{d+1}c_{i,0}X^i=
c_{0,0}+c_{d+1,0}X^{d+1}$ so $a_{1,0}=-\a a_{0,0}$ and
$b_{i,0}=\a^ib_{0,0}$ for $0\leq i\leq d$ for some $\a\in
K^*$. Similarly $(a_{0,1}+a_{1,1}X)\sum_{i=0}^db_{i,1}X^i=
\sum_{i=0}^{d+1}c_{i,2}X^i= c_{0,2}+c_{d+1,2}X^{d+1}$ so $a_{1,1}=-\b
a_{0,1}$ and $b_{i,1}=\b^ib_{0,1}$ for $0\leq i\leq d$ for some $\b\in
K^*$. 

We have $(a_0'+a_1'Y)(b_{0,0}+b_{0,1}Y)=
c_{0,0}'+c_{0,1}'Y+c_{0,2}'Y^2=c_{0,0}'+c_{0,2}'Y^2$ so
$\h{a_1'}{a_0'}=-\h{b_{0,1}}{b_{0,0}}$. Similarly
$(a_0'+a_1'Y)(b_{d,0}+b_{d,1}Y)=c_{d,0}'+c_{d,1}'Y+c_{d,2}'Y^2=
c_{d,0}'+c_{d,2}'Y^2$ so
$\h{a_1'}{a_0'}=-\h{b_{d,1}}{b_{d,0}}=-\h{\b^db_{0,1}}{\a^db_{0,0}}$. So
$(\b/\a )^d=1$. Since $\car K=p$ and $d=p^td'$ this implies $(\b/\a
)^{d'}=1$. It follows that
$\h{b_{d',1}}{b_{d',0}}=\h{\b^{d'}b_{0,1}}{\a^{d'}b_{0,0}}=
\h{b_{0,1}}{b_{0,0}}=-\h{a_1'}{a_0'}$. Since
$(a_0'+a_1'Y)(b_{d',0}+b_{d',1}Y)= c_{d',0}'+c_{d',1}'Y+c_{d',2}'Y^2$
we get $c_{d',1}'=0$. Contradiction. ($c_{i,1}'=0$ holds only for
$i=0,d$.)

Hence $B(I',I'',I,K)$ holds if $\car K\mid d$. We now prove that it
fails otherwise. To do this we have to show that if $\car K\nmid d$
there are $P,Q\in K[X,Y,Z]$ with $I(P)=I'$, $I(Q)=I''$ and
$I(PQ)=I$. We keep the above notation for the coefficients of
$P,Q,PQ$. Since $\car K\nmid d$ we have $|\mu_d(K)|=d$. Let $\zeta$ be
a primitive $d$-th root of unity in $K$ and let $a\in K^*$ with
$a\notin\mu_d\cup\{ -1,-\zeta\1\}$. We take $a_{0,0}=1$, $a_{1,0}=-1$,
$a_{0,1}=a$, $a_{1,1}=-a\zeta$, $a_0'=1$, $a_1'=-1$ and
$b_{i,j}=\zeta^{ij}$ $\forall i,j$. So $P=1-X+a(1-\zeta X)Y+(1-Y)Z$
and $Q=\sum_{i=0}^dX^i+\sum_{i=0}^d\zeta^iX^iY$. Now
$$\sum_{i=0}^{d+1}c_{i,0}X^i=(1-X)\sum_{i=0}^dX^i=1-X^{d+1},$$ 
$$\sum_{i=0}^{d+1}c_{i,2}X^i=a(1-\zeta
X)\sum_{i=0}^d\zeta^iX^i=a(1-\zeta X^{d+1})$$
and
$$\sum_{i=0}^{d+1}c_{i,1}X^i=(1-X)\sum_{i=0}^d\zeta^iX^i+a(1-\zeta
X)\sum_{i=0}^dX^i.$$ 
By the way $a$ was chosen $c_{0,1}=1+a\neq 0$,
$c_{d+1,1}=-1-a\zeta\neq 0$ and for $1\leq i\leq d$ we have
$c_{i,1}=\zeta^i-\zeta^{i-1}+a(1-\zeta)=(\zeta -1)(\zeta^{i-1}-a)\neq
0$. Hence $c_{i,j}=0$ iff $1\leq i\leq d$ and $j=0,2$. 

For $0\leq i\leq d$ we have
$$c_{i,0}'+c_{i,1}'Y+c_{i,2}'Y^2=(1-Y)(1+\zeta^iY)$$ 
so $c_{i,0}=1$, $c_{i,1}'=\zeta^i-1$ and $c_{i,2}'=-\zeta^i$. It
follows that $c_{i,j}'=0$ iff $(i,j)\in\{ (0,1),(d.1)\}$. 

In conclusion, $I(PQ)=I$, as claimed.\qed

{\bf End of proof for Theorem 4.2.} If $S$ is a finite set of primes,
$S\neq\emptyset$ then let $d=\prod_{p\in S}p$. If $I',I'',I\sb\NN^3$
are the sets defined in Lemma 4.4 then $B(I',I'',I,K)$ holds iff $\car
K\mid d$, i.e. iff $\car K\in S$. This implies that if $J\sb\NN^5$ is
defined like in Lemma 4.3 then $A(J,K)$ holds iff $\car K\in S$ so $J$
satisfies the case (a) of Theorem 4.1 for the finite set $S$. \qed
\vskip 3mm

The set $J$ above is included in $\NN^5$,
$|J|=|I|+|I'|+|I''|+1=(4d+7)+6+(2d+2)+1=6d+16$ and $\max\{\Sigma i\mid
i\in J\} =d+3$. So the polinomials in are in $5$ variables and of
degree $d+3$ and have $6d+16$ monomials. We will show that one can
reduce the number of variables to $2$, at the expense of increasing
the degree to $6d+15$. 

\blm Let $I',I'',I\sb\NN^n$ be finite nonempty with $I\sbq I'+I''$. Let
$\ell :\ZZ^n\to\ZZ^m$ be a linear function such that $\ell_{|I'+I''}$
is injective and let $a',a''\in\ZZ^m$ such that $a'+\ell (I'),a''+\ell
(I'')\sb\NN^m$. Then $B(I',I'',I,K)$ is equivalent to $B(a'+\ell
(I'),a''+\ell (I''),a'+a''+\ell (I),K)$. 
\elm
\pf Let $Y=(X_1,\ldots,X_m)$ so $K[Y]=K[X_1,\ldots,X_m]$. 

If $i,i'\in I'$ with $\ell (i)=\ell (i')$ then for any $j\in I''$
$\ell (i+j)=\ell (i'+j)$. Since $i+j,i'+j\in I'+I''$ this implies that
$i+j=i'+j$ so $i=i'$. Thus $\ell_{|I'}$ is injective and similarly for
$\ell_{|I''}$. Therefore we have the bijections $f':W_{I'}\to
W_{a'+\ell (I')}$, $f'':W_{I''}\to W_{a''+\ell (I'')}$ and $f:W_{I'+I''}\to
W_{a'+a''+\ell (I'+I'')}$, given by $\sum_{i\in
I'}a_iX^i\mapsto\sum_{i\in I'}a_iY^{a'+\ell (i)}$, $\sum_{j\in
I''}b_jX^j\mapsto\sum_{j\in I''}b_jY^{a''+\ell (j)}$ $\sum_{k\in
I'+I''}c_kX^k\mapsto\sum_{k\in I'+I''}c_kY^{a'+a''+\ell (k)}$. We have
$f'(V_{I'})=V_{a'+\ell (I')}$, $f''(V_{I''})=V_{a''+\ell (I'')}$ and
$f(V_I)=V_{a'+a''+\ell (I)}$. 

If $P\in W_{I'}$, $Q\in W_{I''}$ then $PQ\in W_{I'+I''}$. If
$P=\sum_{i\in I'}a_iX^i$, $Q=\sum_{j\in I''}a_jX^j$ and $PQ=\sum_{k\in
  I'+I''}c_kX^k$ then $c_k=\sum_{i,j}a_ib_j$, where the sum is taken
over all $i\in I'$, $j\in I''$ with $i+j=k$. 

Now $f'(P)\in W_{a'+\ell (I')}$, $f''(Q)\in W_{a''+\ell (I'')}$ so
$f'(P)f''(Q)\in W_{a'+a''+\ell (I'+I'')}$ so $f'(P)f''(Q)=\sum_{k\in
I'+I''}c_k'Y^{a'+a''+\ell (k)}$ for some $c_k'\in K$. This implies
that for any $k\in I'+I''$ $c_k'=\sum_{i,j}a_ib_j$, where the sum is
taken over all $i\in I'$, $j\in I''$ with $(a'+\ell (i))+(a''+\ell
(j))=a'+a''+\ell (K)$. Since $i+j,k\in I'+I''$ this is equivalent to
$i+j=k$. Hence $c_k'=c_k$ so $f'(P)f''(Q)=f(PQ)$. 

Suppose now that $B(I',I'',I,K)$ fails so there are $P\in V_{I'}$,
$Q\in V_{I''}$ such that $PQ\in V_I$. Then $f'(P)\in V_{a'+\ell
(I')}$, $f''(Q)\in V_{a''+\ell (I'')}$ and $f'(P)f''(Q)=f(PQ)\in
V_{a'+a''+\ell (I)}$ so  $B(a'+\ell (I'),a''+\ell (I''),a'+a''+\ell
(I),K)$ fails. 

Conversely if $B(a'+\ell (I'),a''+\ell (I''),a'+a''+\ell (I),K)$ fails
then let $\( P\in V_{a'+\ell (I')}$, $\( Q\in V_{a''+\ell (I'')}$ such
that $\( P\,\( Q\in V_{a'+a''+\ell (I)}$. Then $\( P=f'(P)$ and $\(
Q=f''(Q)$ for some $P\in V_{I'}$, $Q\in V_{I''}$. Since
$f(PQ)=f'(P)f''(Q)=\( P\,\( Q\in V_{a'+a''+\ell (I)}$ we have $PQ\in
V_I$ so $B(I',I'',I,K)$ fails. \qed

\blm If $I=\prod_{t=0}^s\{ 0,\ldots,b_t-1\}$, where $b_t\geq 2$ are
integers and $\ell :\ZZ^{s+1}\to\ZZ$ is given by
$(i_0,\ldots,i_s)\mapsto i_0+b_0i_1+b_0b_1i_2+\cdots +b_1\cdots
b_{s-1}i_s$ then $\ell$ defines a bijection between $I$ and $\{
0,\ldots,b_0\cdots b_s-1\}$. In particular, $\ell_{|I}$ is injective.
\elm
\pf We use induction on $s$. For $s=0$ our statement is
trivial. Let now $s\geq 1$. Then $I=I'\ti\{ 0,\ldots,b_s-1\}$, where
$I'=\prod_{t=0}^{s-1}\{ 0,\ldots,b_t-1\}$. By the induction hypothesis
the linear mapping $\ell':\ZZ^s\to\ZZ$ given by
$(i_0,\ldots,i_{s-1})\mapsto i_0+b_0i_1+b_0b_1i_2+\cdots +b_1\cdots
b_{s-2}i_{s-1}$ defines a bijection between $I'$ and $\{
0,\ldots,b_0\cdots b_{s-1}-1\}$. 

Now $I=\coprod_{j=0}^{b_s-1}I'\ti\{ j\}$. But for any
$i=(i',j)\in I'\ti\{ j\}$ we have $\ell (i)=\ell'(i')+b_0\cdots
b_{s-1}j$ so $\ell$ defines a bijection between $I'\ti\{ j\}$ and $\{
0,\ldots,b_0\cdots b_{s-1}-1\} +b_0\cdots b_{s-1}j=\{ b_0\cdots
b_{s-1}j,\ldots,b_0\cdots b_{s-1}(j+1)-1\}$. Since
$\coprod_{j=0}^{b_s-1}\{ b_0\cdots b_{s-1}j,\ldots,b_0\cdots
b_{s-1}(j+1)-1\}=\{ 0.\ldots,b_0\cdots b_s-1\}$ we get the desired
result. \qed

Note that if $b_0=\ldots =b_s=b$ then Lemma 4.6 amounts to writing
numbers $\leq b^{s+1}-1$ in base $b$. Also note that we have similar
results if we permutate the roles of the indices $0,\ldots,s$. 

\bff We now consider the sets $I',I'',I\sb\ZZ^3$ from Lemma 4.4. Since
$I'+I''\sb\{ 0,\ldots,d+1\}\ti\{ 0,1,2\}\ti\{ 0,1\}$ by Lemma
4.6 $\ell :\ZZ^3\to\ZZ$ given by $(i,j,k)\mapsto 6i+j+3k$ is injective
on $I'+I''$. Since $\ell (I'),\ell (I'')\sb\NN$ we may apply Lemma
4.5 with $a'=a''=0$. Hence the result from Lemma 4.4 remain true if we
replace $I',I''$ and $I$ by $\ell (I'),\ell (I'')$ and $\ell (I)$. 

That is, if
$$I'=\{ 0,1,3,4,6,7\},$$
$$I''=\{ 6i\mid 0\leq i\leq d\}\cup\{ 6i+1\mid 0\leq i\leq d\},$$
$$I=\{ 0,2,6d+6,6d+8\}\cup\{ 6i+1\mid 0\leq i\leq d+1\}
\cup\{ 6i+3\mid 0\leq i\leq d\}$$
$$\cup\{ 6i+4\mid 1\leq i\leq d-1\}\cup\{
6i+5\mid 0\leq i\leq d\}$$
then $B(I',I'',I,K)$ holds iff $\car K\mid d$. 

We have $\min I'=\min I''=\min I=0$, $\max I'=7$, $\max I''=6d+1$ and
$\max I=6d+8$.
\eff

Since $I',I'',I\sb\NN$ by using Lemma 4.3 one can obtain a set
$J\sb\NN^3$ with $V_J(K)=\emptyset$ iff $\car K\mid d$. The following
lemma proves that in fact we can take $J\sb\NN^2$. 

\blm Let $I',I'',I\sb\NN$ with $\min I'=\min I''=\min I=0$, $\max
I'=g'>0$, $\max I''=g''>0$ and $\max I=g'+g''$. Let $K$ be an
algebrically closed field.

(i) If $h>g'$, $J_0=h+I$, $J_1=I'\cup (h+I'')$ and $J=(J_0\ti\{
0\} )\cup (J_1\ti\{ 1\} )\cup\{(0,2)\}$ then $B(I',I'',I,K)$ is
equivalent to $A(J,K)$. 

(ii) If $I\cup\{ g'\}\neq I'\cup (g'+I'')$ $J_0=g'+I$, $J_1=I'\cup
(g'+I'')$ or $I'\triangle (g'+I'')$  and $J=(J_0\ti\{ 0\} )\cup
(J_1\ti\{ 1\} )\cup\{(0,2)\}$ then $B(I',I'',I,K)$ is equivalent to
$A(J,K)$. 
\elm
\pf We denote $(X_1,X_2)=(X,Y)$. 

We prove (i) and (ii) together. For (ii) we make the convention that
$h=g'$ so in both cases $h\geq g'$. We have $\min I'=0$, $\max I'=g'$,
$\min (h+I'')=h$ and $\max (h+I'')=h+g''$. In the case of (i), when
$g'<h$, we have $I'\cap (h+I'')=\emptyset$, $\min J_1=0$, $\max
J_1=g'+h$. Also $I'=J_1\cap [0,g']$ and $h+I''=J_1\cap [h,h+g'']$. In
the case of (ii), when $h=g'$, we have $I'\cap (g'+I'')=\{ g'\}$ so
$I'\triangle (g'+I'')=I'\cup (g'+I'')\setminus\{ g'\}$. The relations
$\min J_1=0$, $\max J_1=g'+h=g'+g''$ are preserved, even when
$J_1=I'\cup (g'+I'')\setminus\{ g'\}$, since $g'\neq 0,g'+g''$. But
this time $I'=(J_1\cap [0,g'-1])\cup\{ g'\}$ and $g'+I''=\{g'\}\cup
(J_1\cap [g'+1,g'+g''])$. 

Assume that $B(I',I'',I,K)$ fails so there are $P,Q\in K[X]$ with
$I(P)=I'$, $I(Q)=I''$ and $I(PQ)=I$. We take
$F=(\l P+Y)(X^hQ+Y)=\l X^hPQ+(\l P+X^hQ)Y+Y^2$ for some $\l\in
K^*$. Then
$$I(F)=(I(\l X^hPQ)\ti\{ 0\})\cup (I(\l P+X^hQ)\ti\{ 1\})\cup\{
(0,2)\}.$$
We have $I(\l X^hPQ)=h+I=J_0$. In the case (i) since $I(\l P)=I'$ and
$I(X^hQ)=h+I''$ are disjoint we have $I(\l P+X^hQ)=I'\cup
(h+I'')=J_1$. In the case of (ii) we have $I(\l P)\cap I(X^hQ)=I'\cap
(g'+I'')=\{ g'\}$ so $I(\l P+X^hQ)$ equals $I'\cup (g'+I'')\setminus\{
g'\} =I'\triangle (g'+I'')$ if the coefficents of $X^{g'}$ in $\l P$
and $X^{g'}Q$ cancel each other and it equals $I'\cup (g'+I'')$
otherwise. For suitable choices of $\l\in K^*$ each of the two cases
can occur. In particular, $\l$ can be chosen such that $I(\l
P+X^{g'}Q)=J_1$. In conclusion, $I(F)=J$ and since $F$ is decomposable
we have $Z_J\neq\emptyset$ so $A(J,K)$ fails. 

Conversely, assume that $A(J,K)$ fails so there is a decomposable
$F\in K[X,Y]$ with $I(F)=J$. Then $F=R_0+R_1Y+Y^2$, where $R_0,R_1\in
K[X]$ with $I(R_0)=J_0$, $I(R_1)=J_1$. Since $F$ is decomposable we
have $F=(P+Y)(Q+Y)$ so $R_0=PQ$, $R_1=P+Q$ for some $P,Q\in K[X]$. We
have $I(R_0)=h+I$ so $R_0=X^hR_o'$, where $I(R_0')=I$ so $X\nmid R_0'$
and $\deg R_0'=g'+g''$. Since $I(R_1)=J_1$, $\min J_1=0$, $\max
J_1=h+g''$ we have $X\nmid R_1$ and $\deg J_1=h+g''$. Now $X\nmid
R_1=P+Q$ so we may assume that $X\nmid P$. Since $PQ=R_0=X^hR_0'$ we
have $Q=X^hQ'$ for some $Q'\in K[X]$  and $PQ'=R_0'$. We will prove
that $I(P)=I'$ and $I(Q)=h+I''$ so $I(Q')=I''$. Since also
$I(PQ')=I(R_0')=I$ we have that $B(I',I'',I,K)$ fails. 

We prove that $\deg P<h+g''$. In the case of (i) this follows from
$\deg P\leq R_0'=g'+g''$. In the case of (ii) assume that $\deg
P=g'+g''=h+g''$. Then $PQ'=R_0'$ implies that $Q'=\l$ so $Q=\l X^{g'}$
for some $\l\in K^*$. Since $\l P=R_0'$ we have
$I(P)=I(R_0')=I$. Since also $R_1=P+Q=P+\l X^{g'}$ and $I(R_1)=J_1$ we
get $I\cup\{ g'\}=J_1\cup\{ g'\} =I'\cup (g'+I'')$, which contradicts
the hypothesis. So in all cases $\deg P<h+g''=\deg R_0=\deg (P+Q)$. It
follows that $\deg Q=h+g''$ so $\deg Q'=g''$, which, together with
$\deg PQ'=\deg R_0'=g'+g''$, implies $\deg P=g'$. Thus $\max I(P)=g'$,
$\max I(Q')=g''$ and since $X\nmid R_0'=PQ'$ we have $\min I(P)=\min
I(Q')=0$. Since $Q=X^hQ'$ we have $\min I(Q)=h$, $\max I(Q)=h+g''$. 

We have $I(P+Q)=I(R_1)=J_1$. In the case of (i) we have $\min I(P)=0$,
$\max I(P)=g'<h=\min I(Q)$ and $\max I(Q)=h+g''$ so $I(P)=J_1\cap
[0,g']$, $I(Q)=J_1\cap [h,h+g'']$. In the case of (ii) $\min I(P)=0$,
$\max I(P)=g'=\min I(Q)$ and $\max I(Q)=g'+g''$ so $I(P)=(J_1\cap
[0,g'-1])\cup\{ g'\}$, $I(Q)=\{ g'\}\cup (J_1\cap [g'+1,g'+g''])$. In
both cases $I(P)=I'$, $I(Q)=h+I''$, as claimed. \qed

If $I',I'',I\sb\NN$ are the sets defined in 4.7 we have $\min I'=\min
I''=\min I=0$, $\max I'=7$, $\max I''=6d+1$ and $\max I=6d+8$ so we
may apply Lemma 4.8 with $g'=7$, $g''=6d+1$. Since $2\in I$ but
$2\notin I'\cup (7+I'')$ we have $I\cup\{ 7\}\neq I'\cup (7+I'')$ so
we may apply Lemma 4.8(ii). We take $J_0=7+I$, $J_1= I'\triangle
(7+I'')=I'\cup (7+I'')\setminus\{ 7\}$ and $J=(J_0\ti\{ 0\} )\cup
(J_1\ti\{ 1\} )\cup\{ (0,2)\}$. By Lemma 4.8(ii) and 4.7 we have that
$A(J,K)$ iff $B(I',I'',I,K)$ iff $\car K\mid d$. 

Note that $\max J_0=7+\max I=6d+15$ and $\max
J_1=g'+g''=7+(6d+1)=6d+8$ so all polynomials in $V_J$ have degree
$6d+15$. Also $|J_0|=|I|=4d+7$ and
$|J_1|=|I'|+|7+I''|-2=6+(2d+2)-2=2d+6$ so $|J|=(4d+7)+(2d+6)+1=6d+14$
so all polynomials in $V_J$ are sums of $6d+14$ monomials.

\section*{References}

Gao, Shuhong, {\it Absolute irreductibility of polynomials
via Newton polytopes}, J. Algebra 237, No.2, 501-520 (2001).

Kesh, Deepanjan \& Mehta, S. K., {\it Polynomial Irreducibility Testing
through Minkowski Summand Computation},  20th Canadian Conference on
Computational Geometry (CCCG'08), McGill University, Montreal, Canada,
13-15 August 2008. 
\vskip 1cm
\begin{center}
Institute of Mathematics ``Simion Stoilow'' of the Romanian Academy

P.O. Box 1-764, RO-70700 Bucharest, Romania

Email: Constantin.Beli@imar.ro raspopitu1@yahoo.com\end{center}

\end{document}